\documentclass[a4paper, 12pt]{article}

\usepackage[koi8-r,cp1251]{inputenc}
\usepackage{amsfonts,amssymb,amsmath, hyperref}
\usepackage[final]{epsfig}
\usepackage{graphicx}

\begin{document}

\begin{center}
\bf{{\Large On  Hausdorff  operators
on homogeneous   spaces of locally compact groups }}\\
\end{center}

\begin{center}
A. R. Mirotin
\end{center}

\begin{center}
amirotin@yandex.ru
\end{center}

Abstract. Hausdorff  operators on the real line and multidimensional Euclidean spaces originated from some classical summation methods.  Now it is an active research area. Hausdorff  operators on general groups  were defined and studied by the author  since 2019. The purpose of this paper is to define and study Hausdorff  operators on Lebesgue and real Hardy  spaces over homogeneous   spaces of locally compact groups. We introduce in particular an atomic Hardy space over homogeneous   spaces of locally compact groups  and obtain conditions for  boundedness of Hausdorff  operators on such spaces. Several corollaries  are considered and unsolved problems are formulated.

Key wards. Hausdorff  operator, locally compact group, homogeneous   space, atomic Hardy space, Lebesgue space, bounded operator.

\section{Introduction and preliminaries}

 Hausdorff  operators  were introduced by Hardy \cite[Chapter XI]{H} on the segment,
 and by by C.~ Georgakis \cite{G} and independently by E.~ Liflyand and F.~ Moricz \cite{LM} on the whole real line.
 Their multidimensional generalizations were considered later by  Liflyand and Lerner \cite{LL}.
 Now it is an active research area. It is enough to note that the Google search by request "Hausdorff  operator" gives more then  1 200 000 results. See also survey articles \cite{Ls}, \cite{CFW}  for historical remarks and the state of the art up to 2014.

 Hausdorff  operators on general groups  were defined and studied by the author in \cite{JMAA} and \cite{AddJMAA}.
 The purpose of this paper is to define and study Hausdorff  operators on Lebesgue and real Hardy  spaces over homogeneous   spaces of locally compact groups.

 In what follows $G$ stands for a locally compact group with left Haar  measure $\nu$.
 We denote by $\mathrm{Aut}(G)$  the space  of all topological  automorphisms of $G$ endowed with its natural topology (see, e.g. \cite{HiR}),  $\mathcal{L}(Y)$ denotes the space of linear bounded operators on a normed space $Y$.

 In \cite{JMAA} the next definition was proposed.

 \textbf{Definition} \cite{JMAA}. Let  $(T,m)$ be  a measure space, $G$ a topological group, $A:T\to \mathrm{Aut}(G)$ a measurable map,
   and $\Psi$ a locally integrable function on $T.$
   We define the \textit{Hausdorff  operator} with the kernel $\Psi$ over the group $G$  by the formula
$$
\mathcal{H}f(x)=\int_T \Psi(t)f(A(t)(x))dm(t).
$$

By \cite[Lemma 1]{JMAA} this operator is bounded in $L^p(G)$ ($1\le p\le\infty$) provided $\Psi(t)(\mathrm{mod} A(t))^{-1/p}\in L^1(T,m)$ and
$$
\|\mathcal{H}\|_{\mathcal{L}(L^p(G))}\le \int_T |\Psi(t)|(\mathrm{mod} A(t))^{-1/p}dm(t). 
$$

 Moreover, in \cite{JMAA}, \cite{AddJMAA} conditions for boundedness of Hausdorff  operators on the real Hardy space $H^1(G)$ over metrizable locally compact group $G$ with doubling condition were obtained.

 The aim of this work is to define Hausdorff  operators on homogeneous   spaces of locally compact groups and to prove analogs of aforementioned results for this situation.

   Let $K$ be a compact subgroup of  $G$ with normalized Haar measure $\beta$.  Consider the quotient space
$G/K$ of left cosets $\dot x:=xK=\pi_K(x)$ ($x\in G$) where $\pi_K:G\to G/K$ stands for a natural projection. We shall assume that the  measure $\nu$ is normalized in such a way that (generalized) Weil's formula
$$
\int_Gg(x)dx=\int_{G/K}\left(\int_Kg(xk)dk\right)d\lambda(\dot x) \eqno(1)
$$
holds for all $g\in L^1(G)$, where $\lambda$ denotes some left-$G$ invariant measure on $G/K$ (see, e.~g.,
\cite[Chapter VII, §2, No. 5, Theorem 2 ]{Burb} and especially  remark c) after this theorem). \footnote{$G$-left invariance of   $\lambda$ means that $\lambda(xE)=\lambda(E)$ for every Borel subset $E$ of $G/K$ and for every $x\in G$. This measure is unique up to constant multiplier.}
Here and below  we  write $dx$ instead of $d\nu(x)$ and $dk$ instead of $d\beta(k)$. We shall write also $d\dot x$ instead of $d\lambda(\dot x)$.

The function $g:G\to \mathbb{C}$ is called \textit{right-$K$-invariant}  if $g(xk)=g(x)$ for all $x\in G$, $k\in K$. For such a function we put $\dot g(\dot x):=g(x)$. This definition is correct and for $g\in L^1(G)$ formula (1) implies that
$$
\int_Gg(x)dx=\int_{G/K}\dot g(\dot x)d\dot x \eqno(2)
$$
 (recall that $\int_Kdk=1$).

 The map $g\mapsto\dot g$ is a bijection between the set of all right-$K$-invariant functions on $G$ (all right-$K$-invariant functions from $L^1(G)$) and the set of all functions on $G/K$ (respectively functions from $L^1(G/K,\lambda)$).

Let an automorphism $A\in \mathrm{Aut}(G)$ maps $K$ onto itself. Since
$$
A(\dot x):= A(xK)=\{A(x)A(k): k\in K\}=
A(x)K=\pi_K(A(x))
$$
 we get a homeomorphism   $\dot A:G/K\to G/K,$  $\dot A(\dot x):=\pi_K(A(x)).$ Then for every
 right-$K$-invariant function $g$ on $G$ we have  $\dot g(\dot A(\dot x))=g(A(x)).$

From now on we put
 $$
 \mathrm{Aut}_K(G):=\{\dot A : A\in \mathrm{Aut}(G), A(K)=K\}.
 $$

\textbf{Definition 1.} Let  $(\Omega,\mu)$ be  a measure space,  $(\dot A(u))_{u\in \Omega}\subset \mathrm{Aut}_K(G)$  a family of homeomorphisms of $G/K$,  and
$\Phi\in L^1_{loc}(\Omega,\mu).$  For   a  function $f$ on $G/K$ we define a  \textit{Hausdorff operator on} $G/K$ as follows
$$
(\mathcal{H}_{\Phi,\dot A} f)(\dot x):=\int_{\Omega} \Phi(u)f(\dot A(u)(\dot x))d\mu(u).
$$

As was mentioned by Hardy in the case   $\Omega=[0,1]$  \cite[Theorem 217]{H} his Hausdorff  operator possesses some regularity property. A Hausdorff operator in the sense of  Definition 1 also enjoys    this property as the next proposition shows.

\textbf{Proposition 1.} \textit{Suppose that the conditions of Definition 1 are fulfilled and the group $G$ is $\sigma$-compact.
In order that the transformation $\mathcal{H}_{\Phi,\dot A}$  should be
regular, i.e. that $f\in C(G/K),$ $f(\dot x) \to l$ when $\dot x\to \infty$   should imply $\mathcal{H}_{\Phi, \dot A}f(\dot x) \to l$,
 it is necessary and
sufficient that $\int_{\Omega} \Phi(u)d\mu(u)=1.$
}

Proof. If $f(x)=1$ then $\mathcal{H}_{\Phi, \dot A}f(\dot x) =\int_\Omega \Phi(u)d\mu(u).$ Thus, $\int_\Omega \Phi(u)d\mu(u)=1$ is a necessary condition.

To prove the sufficiency, first note that since $\dot A(u)$ has continuous  inverse, $f(\dot x) \to l$ when $\dot x\to \infty$   implies $f(\dot A(u)\dot x)\to l$ when $\dot x\to \infty.$  Indeed,  $f(\dot x) \to l$ when $\dot x\to \infty$ means that  for every $\varepsilon>0$ there is such compact $C_\varepsilon\subset G/K$ that $|f(\dot x)-l|<\varepsilon$ for $x\in G\setminus C_\varepsilon$. Now if $\dot x\in (G/K)\setminus \dot A(u)^{-1}(C_\varepsilon)$ we get  $\dot A(u)\dot x\in (G/K)\setminus C_\varepsilon$ and therefore $|f(\dot A(u)\dot x)-l|<\varepsilon$.

But if, in addition,  $f\in C(G/K)$ the function $f$ is bounded on $G/K$ and therefore  $\mathcal{H}_{\Phi,\dot A}f(\dot x) \to l$  by the Lebesgue Theorem (one can apply the Lebesgue Theorem, since $G/K$ is $\sigma$-compact).

Thus, Proposition 1 shows that Hausdorff operators in a sense of Definition 1 gives us a family (for various $\Phi$, $A(u)$, and $\Omega$) of generalized limits at infinity for functions on $G/K$.

\section{Hausdorff  operators on $L^p(G/K)$}

In the following we put $L^p(G/K):=L^p(G/K,\lambda)$ ($p\in [1,\infty]$). Formula (2) implies that $\|g\|_{L^p(G)}=\|\dot g\|_{L^p(G/K)}$. For every right-$K$-invariant function $g\in L^p(G)$.

In this section we give  conditions of boundedness of Hausdorff operators on $L^p(G/K)$.

Let $(\Omega,\mu)$ and $(\dot A(u))_{u\in \Omega}$ be as in Definition 1. For a function $\Phi$ on $\Omega$ let
$$
\|\Phi\|_{p,A}:=\int_{\Omega}|\Phi(u)|(\mathrm{mod} A(u))^{-1/p}d\mu(u).
$$

\textbf{Theorem 1.} \textit{Suppose that the conditions of Definition 1 are fulfilled, $p\in [1,\infty]$, and $\|\Phi\|_{p,A}<\infty$. Then $\mathcal{H}_{\Phi, \dot A}$ is bounded in $L^p(G/K)$ and}
$$
\|\mathcal{H}_{\Phi, \dot A}\|_{\mathcal{L}(L^p(G/K))}\le \|\Phi\|_{p,A}.
$$

Proof. Let $1< p<\infty$. Every function $f\in L^p(G/K)$ has the form $f=\dot g$ for a unique right $K$-invariant function $g\in L^p(G)$. Using Minkowski’s integral inequality we have   that
$$
\|\mathcal{H}_{\Phi, \dot A}f\|_{L^p(G/K)}=\left(\int_{G/K} \left|\int_{\Omega}\Phi(u)f(\dot A(u)(\dot x))d\mu(u) \right|^pd\dot x\right)^{1/p}
$$
$$
\leq\int_{\Omega}|\Phi(u)|\left(\int_{G/K}|f(\dot A(u)(\dot x))|^pd\dot x\right)^{1/p}d\mu(u).
$$
Since the function $x\mapsto g(A(u)(x))$ is right $K$-invariant, as well, formula (2) yields
$$
\int_{G/K}|f(\dot A(u)(\dot x))|^pd\dot x=\int_{G}|g(A(u)(x))|^pdx.
$$
On the other hand, by \cite[VII.1.4, formula (31)]{Burb} we have
$$
\int_G|g(A(u)(x))|^pdx=(\mathrm{mod} A(u))^{-1}\int_G|g(x)|^pdx.
$$
Thus,
$$
\|\mathcal{H}_{\Phi, \dot A}f\|_{L^p(G/K)}\leq\int_{\Omega}|\Phi(u)|(\mathrm{mod} A(u))^{-1/p}d\mu(u)\left(\int_G|g(x)|^pdx\right)^{1/p}
$$
$$
=\|\Phi\|_{p,A}\left(\int_{G/K}|f(\dot x)|^pd\dot x\right)^{1/p}
=\|\Phi\|_{p,A}\|f\|_{L^p(G/K)}.
$$
For $p=1$ the statement of  theorem 1 follows from Fubini theorem and for $p=\infty$ it is obvious.

The following simple example is intended to illustrate preceding  constructions.

\textbf{Example 1.} Let $G$ be the multiplicative group $\mathbb{C}^\times$ of the complex field $\mathbb{C}$ and $K:=\{z\in \mathbb{C}^\times: |z|=1\}$ the circle group. Then $G/K$ can  be identified with $(0,\infty)$ via the map $zK= \dot z\mapsto r$ where $r=|z|$. In other words, we use $(0,\infty)$ as a model of   $G/K$, the positive number $r$ representing the circle of radius $r$.
Automorphisms of $G$ have the form $A(re^{i\alpha})= r^ue^{i\alpha}$ or $A(re^{i\alpha})= r^ue^{-i\alpha}$ ($u\in \mathbb{R}$). Thus $\mathrm{Aut}(G)/K=\mathrm{Aut}(G)$.
It follows that $\dot A(\dot z)=r^u$  ($u\in \mathbb{R}$).
 So  the general form of a Hausdorff operator on $(0,\infty)$ looks as follows ($f:(0,\infty)\to \mathbb{C}$, $r>0$)
$$
\mathcal{H}_{\Phi, \dot A}f(r)=\int_{\mathbb{R}}\Phi(u)f(r^u)d\mu(u)
$$
 (we take $\Omega=\mathbb{R}$, and $\mu$  any positive measure on $\mathbb{R}$). Since $G$ is commutative, $\mathrm{mod} A=1$ for all $A\in \mathrm{Aut}(G)$. So Theorem 1 implies that $\mathcal{H}_{\Phi, \dot A}$ is bounded on $L^p(0,\infty)$ if $\Phi\in L^1(\mu)$ and
in this case $\|\mathcal{H}_{\Phi, \dot A}\|_{\mathcal{L}(L^p)}\le \|\Phi\|_{L^1(\mu)}$.

If we take in Theorem 1 the space $\mathbb{Z}_+$ (endowed  with counting measure) as $\Omega$, we arrive  at the following

\textbf{Corollary 1.} \textit{Let $\Phi(j)$ be a sequence of complex numbers. Consider a discrete Hausdorff  operator over $G/K$}   \footnote{ Discrete Hausdorff  operators were introduced in \cite{Forum}, \cite{faa}.}
$$
\mathcal{H}_{\Phi,\dot A}f(x):=\sum_{j=0}^\infty \Phi(j)f(\dot A(j)(x)).
$$
\textit{Then}
$$
\|\mathcal{H}_{\phi,\dot A}\|_{\mathcal{L}(L^p(G/K))}\le\sum_{j=0}^\infty |\Phi(j)|(\mathrm{mod}A(j))^{-1/p}.
$$

Putting
$$
\Phi=\frac{\chi_{\{\mathrm{mod}A(u)\ge 1\}}}{\mathrm{mod}A}
$$
in Definition 1 ($\chi_E$ denotes the indicator function of the set $E$) we get a \textit{\'{C}esaro  operator} over $G/K$ (cf. \cite{JMAA})
$$
\mathcal{C}_{A,\mu}f(x):=\int_{\{\mathrm{mod}A(u)\ge 1\}}\frac{f(\dot A(u)(x))}{\mathrm{mod}A(u)}d\mu(u).
$$

\textbf{Corollary 2.} \textit{For a \'{C}esaro  operator over $G/K$ the following estimate holds}

$$
\|\mathcal{C}_{A,\mu}\|_{\mathcal{L}(L^p(G/K))}\le \int_{\{\mathrm{mod}A(u)\ge 1\}}\frac{d\mu(u)}{(\mathrm{mod}A(u))^{1+1/p}}.
$$

\section{Hausdorff  operators on Hardy space $H^1(G/K)$}

The  goal  of  this section is to introduce the atomic  Hardy space $H^1(G/K)$ and  to obtain conditions for boundedness of Hausdorff operators on this space.

In the rest of the paper as in \cite {AddJMAA} we assume in addition   that $G$ is  metrizable via a left invariant metric $\rho$ and  the following \textit{doubling condition} in a  sense of \cite{CW} holds:

there exists  a constant $C$ such that
$$
\nu(B(x,2r))\leq C \nu(B(x,r))
$$
for each $x\in G$ and $r > 0$.

 Here $B(x,r)$ denotes the ball of radius $r$ around $x.$ The \textit{doubling constant }  is the smallest constant
$C\geq  1$ for which the last inequality holds. We denote this constant by $C_\nu.$ Then  for each $x\in G, k\geq 1$ and $r > 0$
$$
\nu(B(x,kr))\leq C_\nu k^{s} \nu(B(x,r)),\eqno(3)
$$
where $s=\log_2C_\nu$ (see, e.g., \cite[p. 76]{HK}). The number $s$   takes
the role of a "dimension"  for a doubling metric measure space $G$.

To introduce the space $H^1(G/K)$ first recall  that a function
$a$ on $G$ is an ($(1, \infty)$-)\textit{atom} if

(i) the support of $a$ is contained in a ball $B(x,r)$;

(ii) $\|a\|_\infty\leq\frac{1}{\nu(B(x,r))};$

(iii) $\int_G a(x)dx = 0$.

In case $\nu(G)<\infty$ we shall
 assume $\nu(G)=1$. Then the constant function having value  $1$ is also
considered to be an atom.

\textbf{Definition 2.} We define the Hardy space $H^1(G/K)$ as a space of such functions  $f=\dot g$ on  $G/K$ that $g$ admits an atomic decomposition  of the form
$$
g=\sum_{j=1}^\infty \alpha_ja_j,\eqno(4)
$$
where $a_j$ are  right-$K$-invariant atoms and $\sum_{j=1}^\infty |\alpha_j|<\infty$.
In this case,
$$
\|f\|_{H^1(G/K)}:=\inf\sum_{j=1}^\infty |\alpha_j|,
$$
and infimum is taken over all  decompositions above of $g$.

Thus a function $f=\dot g$ from $H^1(G/K)$ admits an atomic decomposition
$f=\sum_{j=1}^\infty \alpha_j\dot a_j$ such that  $\sum_{j=1}^\infty |\alpha_j|<\infty$, and  $\|f\|_{H^1(G/K)}=\|g\|_{H^1(G)}$.

\textbf{Remark 1}. Real Hardy spaces over compact connected (not necessary quasi-metric) Abelian groups  were  defined in \cite{Indag}. The case of semisimple Lie groups was considered earlier in \cite{Kaw}.

\textbf{Proposition 2}. \textit{The space $H^1(G/K)$ is Banach.  If for some  $h\in H^1(G)$
the inequality    $\int_K h(k)dk\ne 0$ holds the space  $H^1(G/K)$ is nontrivial.}

Proof.  We prove that   right-$K$-invariant atoms exist. To this end for an atom $a$ on $G$ let's consider the function
$$
a'(x):=c\int_K a(xk)dk.
$$
Then $a'$ is right-$K$-invariant and satisfies (i) for every constant $c>0$. Indeed, if $a$ is supported in a ball $B=B(e,r_B),$  then $a'(x)=0$ for $x\notin BK.$ Since
$$
\rho(e,xk)\le \rho(e,x)+\rho(x,xk)=\rho(e,x)+\rho(e,k)< r_B+\mathrm{dist}(e,K)=:r'_B
$$
for $x\in B$, the  set $BK$ is contained in a  ball $B'=B(e,r'_B)$.  Thus, $a'$  is supported in $B'$.
From now on we assume that $c= \nu(B)/\nu(B')$.
Then (ii) holds  for $a'$ because      $\|a'\|_\infty\le c\|a\|_\infty\le c/\nu(B)= 1/\nu(B')$
  for such $c$.
 The property (iii) for $a'$ follows from
  $$
  \int_G\!a'(x)dx=c\!\!\int_G\!\int_K\!a(xk)dkdx=c\!\!\int_K\! \Delta_G(k)\! \int_G \!a(x)dxdk=0.
  $$

Finally, we shall show that $H^1(G/K)$ is complete. First note that since $\|a\|_{L^1(G)}\le 1$ for each atom $a$, we have $\|g\|_{L^1(G)}\le \|\dot g\|_{H^1(G/K)}$ for each right-$K$-invariant function  $g\in H^1(G)$. Then for a function $f=\dot g$ we have $\|f\|_{L^1(G/K)}\le \|f\|_{H^1(G/K)}$.

Let   a sequence  $f_j\in H^1(G/K)$ be such that $\sum_j\|f_j\|_{H^1(G/K)}<\infty$. It is enough to prove that the series $\sum_j f_j$ converges in $H^1(G/K)$. The sequence $S_n$ of  partial sums of this series is a Cauchy sequence in $L^l(G/K)$ because for $m<n$
$$
\|S_n-S_m\|_{L^l(G/K)}\le \sum_{j=m+1}^n \|f_j\|_{L^1(G/K)} \le\sum_{j=m+1}^n \|f_j\|_{H^1(G/K)}.
$$
 So the series $\sum_j f_j$ converges  in $L^1(G/K)$ to some function $f$. On the other hand, each $f_j$ has an atomic decomposition $f_j= \sum_i\alpha_{ij}\dot a_{ji}$ such that $\sum_i|\alpha_{ij}|<2 \|f_j\|_{H^1(G/K)}$. Then $f$ has an atomic decomposition
 $$
 f=\sum_j\sum_i\alpha_{ij}\dot a_{ji},
 $$
 and
 $$
 \|f\|_{H^1(G/K)}\le \sum_j\sum_i|\alpha_{ij}|<2\sum_j\|f_j\|_{H^1(G/K)}<\infty.
 $$
Thus $f\in H^1(G/K)$. Moreover,
$$
\left\|f-\sum_{j=1}^n f_j\right\|_{H^1(G/K)}\le  \sum_{j=n+1}^\infty\|f_j\|_{H^1(G/K)}\to \mbox{ as } n\to\infty.
$$
This completes the proof.

In the proof of Theorem 2 the next lemmas  play an important role.

\textbf{Lemma 1.}  \cite{JMAA} \textit{ Let $(\Omega,q,\mu)$ be  $\sigma$-compact quasi-metric space
with quasi-metric $q$ and  positive Radon measure $\mu,$  $(X,m)$ be a measure space and $\mathcal{F}(X)$
be some Banach  space of  $m$-measurable functions on $X.$ Assume  that the convergence
of a sequence strongly in  $\mathcal{F}(X)$ yields the convergence of some subsequence
to the same function for $m$-almost all $x\in X.$ Let $F(u,x)$ be a function such that
$F(u,\cdot)\in \mathcal{F}(X)$ for $\mu$-almost all $u\in \Omega$ and the
map $u\mapsto F(u,\cdot):\Omega\to \mathcal{F}(X)$ is Bochner integrable with respect
to $\mu$. Then for $m$-almost all $x\in X$}
$$
\left((B)\int_\Omega F(u,\cdot)d\mu(u)\right)(x)=\int_\Omega F(u,x)d\mu(u).
$$

\textbf{Lemma 2.}  \cite{AddJMAA}  \textit{Every automorphism $A\in \mathrm{Aut}(G)$ is Lipschitz. Moreover, one can choose the Lipschitz constant to be
$$
L_A = \kappa_\rho\mathrm{mod}A,
$$
where the constant $\kappa_\rho$ depends on the metric $\rho$  only}.

Now we are in position to prove the following

\textbf{Theorem 2.} \textit{Under assumptions of Definition 1 let $(\Omega,q,\mu)$ be  $\sigma$-compact quasi-metric space
with quasi-metric $q$ and  positive Radon measure $\mu$ and
$\Phi\in L^1(k^s\mu)$ where $k(u):=\kappa_\rho/\mathrm{mod}(A(u))$. Then the operator $\mathcal{H}_{\Phi,\dot A}$ is bounded on the  space $H^1(G/K)$ and }
$$
\|\mathcal{H}_{\Phi,\dot A}\|_{\mathcal{L}(H^1(G/K))}\leq C_\nu\|\Phi\|_{L^1(k^s\mu)}.
$$

Proof. We proceed as in the proof of main theorem in \cite{JMAA}, \cite{AddJMAA}. If we set $X=G/K$  and $m=\lambda$
 the pair $(X,m)$ satisfies  the conditions of Lemma 1 with $H^1(G/K)$ in place of $\mathcal{F}(X)$. Indeed, let $f_n=\dot g_n\in H^1(G/K)$, $f=\dot g\in H^1(G/K)$, and $\|f_n-f\|_{H^1(G/K)}\to 0$ ($n\to\infty$). Since
$$
\|f_n-f\|_{L^1(G/K)}=\int_{G/K}|\pi_K(g_n- g)|d\lambda
$$
$$
=\int_{G}|g_n(x)-g(x)|dx\leq \|g_n-g\|_{H^1(G)}=\|f_n-f\|_{H^1(G/K)}\to 0,
$$
there is a subsequence of $f_{n}$ that converges to $f$ $\lambda$-a.e.

Then Definition 2 and lemma 1 imply for $f\in H^1(G/K)$ that
$$
\mathcal{H}_{\Phi, \dot A}f= \int_\Omega \Phi(u) f\circ \dot A(u)d\mu(u)
$$
(the Bochner integral).

Therefore (below $f=\dot g$)
$$
\|\mathcal{H}_{\Phi,\dot A}f\|_{H^1(G/K)}\leq \int_\Omega |\Phi(u)|\|f\circ \dot A(u)\|_{H^1(G/K)}d\mu(u)
$$
$$
=\int_\Omega |\Phi(u)|\|g\circ A(u)\|_{H^1(G)}d\mu(u).\eqno(5)
$$
If $g$ has representation (4) then
$$
g\circ A(u)=\sum_{j=1}^\infty \alpha_ja_j\circ A(u).\eqno(6)
$$
We claim that $b_{j,u}:=C_\nu^{-1}k(u)^{-s}a_j\circ A(u)$ is an atom, too. Indeed, Lemma 2 implies that
$$
A(u)^{-1}B(x,r))\subseteq B(x',k(u)r),
 $$
 where $x'=A(u)^{-1}(x)$. If $a_j$ is supported in $B(x_j,r_j)$ then $b_{j,u}$ is supported in $B(x_j',k(u)r_j)$. So the condition (i) holds for  $b_{j,u}$.

 Next, by (3) we have
 $$
 \nu(B(x_j',k(u)r_j))\le C_\nu k(u)^s\nu(B(x_j,r_j)).
 $$
This implies that
$$
\|a_j\circ A(u)\|_\infty=\|a_j\|_\infty\le \frac{1}{\nu(B(x_j,r_j))}\le C_\nu k(u)^s\frac{1}{\nu(B(x_j',k(u)r_j))}.
$$
So, the condition (ii) is also fulfilled  for  $b_{j,u}$.

The validity of (iii) follows from \cite[VII.1.4, formula (31)]{Burb}.

Finally, since formula (6) can be rewritten as
$$
g\circ A(u)=\sum_{j=1}^\infty (C_\nu k(u)^s\alpha_j)b_{_ju},
$$
we get
$$
\|g\circ A(u)\|_{H^1(G)}\le C_\nu k(u)^s\sum_{j=1}^\infty |\alpha_j|,
$$
which in turn implies that
$$
\|g\circ A(u)\|_{H^1(G)}\le C_\nu k(u)^s\|g\|_{H^1(G)}=C_\nu k(u)^s\|f\|_{H^1(G/K)}.
$$
Thus, the statement of the theorem  follows from formula (5).

\textbf{Corollary 3.} \textit{Let the assumptions of  Theorem 2 holds.  Then we have for a discrete Hausdorff  operator over $G/K$ (see Corollary 1)}
$$
\|\mathcal{H}_{\Phi,\dot A}\|_{\mathcal{L}(H^1(G/K))}\le C_\nu\kappa_\rho^s\sum_{j=0}^\infty \frac{|\Phi(j)|}{(\mathrm{mod}A(j))^s}.
$$

\textbf{Corollary 4.} \textit{Let the assumptions of  Theorem 2 holds.  Then we have for a \'{C}esaro  operator over $G/K$ (see Corollary 2)}

$$
\|\mathcal{C}_{A,\mu}\|_{\mathcal{L}(H^1(G/K))}\le C_\nu\kappa^s_\rho\int_{\{\mathrm{mod}A(u)\ge 1\}}\frac{d\mu(u)}{(\mathrm{mod}A(u))^{1+s}}.
$$

(Recall, that for this operator
$$
\Phi=\frac{\chi_{\{\mathrm{mod}A(u)\ge 1\}}}{\mathrm{mod}A}.)
$$

\section*{Concluding remarks}
It would be of interest to apply Theorems 1 and 2 to classical homogeneous   spaces  such as Euclidean plane $\mathbb{R}^2=M(2)/O(2),$ sphere $\mathrm{S}^2=O(3)/O(2)$, non-Euclidean plane $\mathbb{H}=SU(1,1)/SO(2)$ \cite[\S 4]{Helg}, to other Riemannian symmetric spaces etc. Would also be of interest  to generalize Theorems 1 and 2 to the case when the group $K$ is noncompact.

This is a preprint of the paper \cite{JBGU}.

\end{document}